\documentclass[final]{siamltex}
\input epsf.tex

\newtheorem{prop}{Proposition}

\newtheorem{defn}{Definition}
\newcommand{\bth}{\begin{theorem}} 
\newcommand{\eth}{\end{theorem}}

\newcommand{\bpr}{\begin{proof}}
\newcommand{\epr}{\end{proof}}

\title{High Accuracy Method for Integral Equations with Discontinuous
Kernels}
\author{SHEON-YOUNG KANG\footnotemark[1]
\and ISRAEL KOLTRACHT\footnotemark[1]
\and GEORGE RAWITSCHER\footnotemark[2]}

\begin{document}

\maketitle

\renewcommand{\thefootnote}{\fnsymbol{footnote}}
\footnotetext[1]{Department of Mathematics, University of Connecticut,
Storrs, CT 06269}
\footnotetext[2]{Department of Physics, 
University of Connecticut, Storrs, CT 06269 }

\begin{abstract}
A new highly accurate numerical approximation scheme 
based on a Gauss type  
Clenshaw-Curtis Quadrature for Fredholm integral equations 
of the second kind
 $$x(t)+\int^{b}_{a}k(t,s)x(s)ds=y(t)$$
whose kernel $k(t,s)$ 
is either discontinuous or 
not smooth along the main diagonal,
is presented.
This scheme is of spectral accuracy when $k(t,s)$ is infinitely differentiable
away from the diagonal $ t = s,$  \ and is also applicable  when $k(t,s)$ \
is singular along the boundary, and at isolated points on the main diagonal.
The corresponding composite rule is described. Application to 
integro-differential Schroedinger equations with non-local potentials is given.
\end{abstract}



\pagestyle{myheadings}
\thispagestyle{plain}
\markboth{S.-Y. KANG, I. KOLTRACHT AND G. RAWITSCHER}{}
 
\section{Introduction}
Let the integral operator,
$$({\it  K}x)(t)=\int^{b}_{a}k(t,s)x(s)ds, \hspace{5mm} 
a \leq t \leq b,$$ map
${\it C}^{q}_{[a, b]}, \  q>1,$ into itself.
In the present paper, we consider the numerical solution of 
the corresponding Fredholm integral equation of the second kind, 
\begin{equation}
\label{es:s1}
x(t)+\int^{b}_{a}k(t,s)x(s)ds=y(t), \hspace{3mm} y \in C^{q},
\hspace{3mm} a \leq t \leq b.
\end{equation}
When the kernel $k(t,s)$ has a discontinuity either by itself 
or in its partial  derivatives along the main diagonal $t = s$,
one can not expect a high accuracy
 numerical approximation based on
Newton-Cotes or  Gaussian Quadratures, \ 
 ( see e.g. Figure 2 of Section 6 ), since,
except for the Trapezium rule,
 the standard error bounds for these rules are not applicable.  
In the present paper we introduce a high accuracy 
discretization technique based on a
Gauss type Clenshaw-Curtis quadrature for a certain class of such
kernels which we call {\bf semismooth}.
\begin{defn}
A kernel $k(t,s)$ is called $(p_{1},p_{2})-semismooth$, if \\
\[
k(t,s)=
\left \{ \begin{array}{cc}
k_{1}(t,s) & {\rm if} \ \  a \leq s \leq t  \\
k_{2}(t,s) & {\rm if} \ \  t \leq s \leq b , 
\end{array}
\right.
\]
where
$k_{1}(t,s) \in  C^{p_{1}}_{[a,b]\times [a,b]}$  
and \  
$k_{2}(t,s) \in  C^{p_{2}}_{[a,b]\times [a,b]}$ \hspace{2mm} for \ some 
$ p_{1}, p_{2} > 1$. 
\end{defn}

Note that for our purpose each of the auxiliary kernels 
$k_1(t,s)$ and $k_2(t,s)$
must be defined in the whole square $[a,b] \times [a,b]$. 
The convergence of our method is of ${\it{O}}(n^{1-r}),$ \  where
$r=\min \{p_{1},p_{2},q\}.$
When $r=\infty,$ the convergence is
superalgebraic,
or spectral. Even for some singular
kernels, when the obtained error estimates are not applicable,
the method still shows good accuracy on numerical examples 
 due to the clustering of mesh 
 points around singularities.  
The 2-step method of Deferred Approach to the Limit, based on the
Trapezium rule,  ( see  e.g. Baker \cite{Baker}, p 363),  works
well for nonsingular kernels, but it is 
much more time consuming, for comparable accuracy, and is not applicable to
kernels with singularities. 

In a way the present method is an extension of a previously examined situation,
which occurs  when a second order differential equation, such as the
Schroedinger equation with local potentials, is transformed
into an integral equation, \cite{reo2}. The kernel of this integral
equation is generally semismooth, but with the additional
structure of $k_{1,2}(t,s)$ being of low rank. Such kernels
we call {\bf semiseparable}. In this case the algorithms developed
in \cite{reo}, \cite{reo2} are adequate and give fast
and accurate solution. However, for Schroedinger equations with
non-local potentials, the corresponding kernels
may be semismooth but not semiseparable, for which our present
technique is perfectly well suited.
This situation
is examined in more detail in Section 7.
  
In Section 2, we describe the discretization of equation (\ref{es:s1})
which is based 
on the Clenshaw-Curtis quadrature for a smooth kernel $k(t,s)$.
This discretization is different from the usual
Gauss-Chebyshev quadrature, (e.g. Delves-Mohamed, \cite{DelMo}), even
for smooth kernels, and from that of Reichel \cite{Reichel}, based on Chebyshev
polynomial expansions. 
In Section 3 we consider semismooth kernels and show that 
the application of Clenshaw-Curtis quadrature 
results in a linear system of equations whose coefficient
 matrix is defined in terms of
  {\bf Schur}, or {\bf componentwise}, products of given matrices.
The accuracy of approximation is determined by the smoothness of $k_{1}$ \ and
$k_{2}$ \ only, and is not affected by the discontinuity along the diagonal $t = s$.
For smooth kernels this discretization is identical with the one described in
Section 2.
The proposed method works well also for the case when a semi-smooth kernel
has singularities on the boundary of the square $[a, b] \times [a, b]$,
although the accuracy is not spectral anymore.
The success of this method in this case is due to the clustering of
Chebyshev support points
near the boundaries, where the singularities occur.
In Section 4 we describe the corresponding composite rule and in Section 5
we apply it to kernels which have finite number of singularities on 
the main diagonal $t=s$.
In Section 6  we describe numerical experiments and comparisons 
with  relevant existing \
methods for kernels with various discontinuities and
singularities. In Section 7 we  
apply the developed technique to the solution of radial Schroedinger 
integro-differential equations 
with a non-local potential.

\section{Discretization of a Smooth Kernel}

Let $k(t,s)$ be differentiable in $t$ and $s$.
Assume that $k(t_{k},s)x(s)$ as a function of $s$ 
 can be expanded in a finite set of polynomials
, i.e.,
\begin{equation}
\label{es:s2}
k(t_{k},s)x(s)=\sum_{j=0}^{n}a_{kj}T_{j}(s), \hspace{3mm}  -1 \leq s \leq 1,
\end{equation}
where \ \ 
$T_{j}(s)=cos(jarccos(s)), \hspace{3mm} j=0,1,...,n,$
are the Chebyshev polynomials. Without any loss of generality
we assume for now that $a=-1$ \ and $b=1$ in equation (\ref{es:s1}).
Let
\begin{eqnarray*}
F(r)=\int^{r}_{-1}k(t_{k},s)x(s)ds=\sum^{n+1}_{j=0}b_{kj}T_{j}(r).
\end{eqnarray*} 
Clenshaw and Curtis \cite{ClCu} showed that  
$$\left[b_{k0},  b_{k1},...,b_{kn+1}\right]^{T}={\bf{S}}_{L}
\left[a_{k0},a_{k1},...,a_{kn}\right]^{T},$$
where
\begin{eqnarray*}
{\bf{S}}_{L}=
\left[\begin{array}{llllll}
1 & 1 & -1 & 1 & \cdots & (-1)^{n} \\
0 & 1 & 0  & 0 & \cdots & 0 \\
0 & 0 & 1 & 0 & \cdots & 0 \\
\vdots & \vdots & \vdots & \ddots & \ddots & \vdots \\
0 & 0 & 0 & 0 & 1& 0 \\
0 & 0 & \cdots & 0 & 0 & 1 \\
\end{array}
\right]
\left[\begin{array}{llllll}
0 & 0 & 0 & 0 & \cdots & 0 \\
1 & 0 & \frac{-1}{2}  & 0 & \cdots & 0 \\
0 & \frac{1}{4} & 0 & \frac{-1}{4} & \cdots & 0 \\
\vdots & \vdots & \ddots & \ddots & \ddots & \vdots \\
0 & \cdots & 0 & \frac{1}{2(n-1)} & 0 & \frac{-1}{2(n-1)} \\
0 & \cdots & 0 & 0 & \frac{1}{2n} & 0 \\
\end{array}
\right]
\end{eqnarray*}
is the so called left spectral integration matrix. 
Here $\left[\nu\right]^{T}$ denotes the transpose of the column vector $\nu.$ 
Since \ \ $T_{j}(1)=1,$ for $j = 0,1,...,n,$ it follows that
\begin{eqnarray*}
F(1) & = & \int^{1}_{-1}k(t_{k},s)x(s)ds=\sum^{n+1}_{j=0}b_{kj} \\
     & = & [1,...,1]\left[b_{k0},  b_{k1},...,b_{kn+1}\right]^{T} 
      =  [1,...,1]{\bf{S}}_{L}\left[a_{k0},a_{k1},...,a_{kn}\right]^{T}.
\end{eqnarray*}
Let $\tau_{k}, k=0,1,...,n,$ denote the zeros of $T_{n+1}$, viz.,\\
$$\tau_{k}=cos\frac{(2k+1)\pi}{2(n+1)},$$
so that $$T_{j}(\tau_{k})=cos\frac{(2k+1)j\pi}{2(n+1)}, \hspace{5mm}
k,j=0,1,...,n.$$
Substituting $s=\tau_{k}, k=0,1,...,n,$ into (\ref{es:s2}), we obtain that
\begin{eqnarray*}
\left[\begin{array}{c}
\alpha_{k0}\\
\alpha_{k1}\\
\vdots\\
\alpha_{kn}
\end{array}
\right]={\bf{C}}^{-1}
\left[\begin{array}{c}
k(t_{k},\tau_{0})x(\tau_{0})\\
k(t_{k},\tau_{1})x(\tau_{1})\\
\vdots\\
k(t_{k},\tau_{n})x(\tau_{n})
\end{array}
\right],
\end{eqnarray*}
where ${\bf{C}}^{-1}$ is an inverse of  {\it{discrete cosine}} transformation 
\ matrix ${\bf{C}}$  whose elements are
specified by 
$${\bf{C}}_{kj}=T_{j}(\tau_{k}), \hspace{3mm} k,j=0,1,...,n.$$
The matrix ${\bf{C}}$ has orthogonal columns, that is,
${\bf{C}}^{T} {\bf{C}}= diag(n,\frac{n}{2},...,\frac{n}{2}).$ 
Therefore,
${\bf{C}}^{-1}=diag(\frac{1}{n},\frac{2}{n},...,\frac{2}{n}){\bf{C}}^{T}$. 
By choosing $t_{k}$ in (\ref{es:s2}) to be Chebyshev points and 
by substituting 
$t=\tau_{k}$ into (\ref{es:s1}),
we get
\begin{eqnarray*}
y(\tau_{k})  =  
x({\tau}_{k})+[1,...,1]{\bf{S}}_{L}{\bf{C}}^{-1}
diag(k(\tau_{k},\tau_{0}),k(\tau_{k},\tau_{1}),...,k(\tau_{k},\tau_{n}))
\left[\begin{array}{c}
x(\tau_{0})\\
x(\tau_{1})\\
\vdots\\
x(\tau_{n})
\end{array}
\right].
\end{eqnarray*}
Introducing \ \ 
$[\sigma_{0}, \sigma_{1},...,\sigma_{n}]=[1,1,...,1]{\bf{S}}_{L}{\bf{C}}^{-1}$
\ we can write
\begin{eqnarray*}
y(\tau_{k}) =
[\sigma_{0}, \sigma_{1},...,\sigma_{n}]
diag(k(\tau_{k},\tau_{0}),k(\tau_{k},\tau_{1}),...,k(\tau_{k},\tau_{n}))
\left[\begin{array}{c}
x(\tau_{0})\\
x(\tau_{1})\\
\vdots\\
x(\tau_{n})
\end{array}
\right],
\end{eqnarray*} 
or equivalently,
\begin{eqnarray*}
y(\tau_{k}) =
[k(\tau_{k},\tau_{0}),k(\tau_{k},\tau_{1}),...,k(\tau_{k},\tau_{n})]
diag(\sigma_{0}, \sigma_{1},...,\sigma_{n})
\left[\begin{array}{c}
x(\tau_{0})\\
x(\tau_{1})\\
\vdots\\
x(\tau_{n})
\end{array}
\right].
\end{eqnarray*} 
Therefore the discretization of the equation (\ref{es:s1}) for the case
$a=-1$ \ and \
$b=1$ 
\ is as follows,
\begin{equation}
\label{es:s3}
\left[{\bf{I}}+{\bf{K}}{\bf{D}}_{\sigma}\right]{\bar{\bf{x}}}={\bar{\bf{y}}},
\end{equation}
where
\begin{eqnarray*}
{\bf{K}} & = & (k(\tau_{i},\tau_{j}))^{n}_{i,j=0} , \\
{\bf{D}}_{\sigma} & = & diag(\sigma_{0}, \sigma_{1},...,\sigma_{n}),\\
{\bar{\bf{x}}} & = &
\left[x(\tau_{0}), x(\tau_{1}),...,x(\tau_{n})\right]^{T},\\
{\bar{\bf{y}}} & = &
\left[y(\tau_{0}), y(\tau_{1}),...,y(\tau_{n})\right]^{T}.\\
\end{eqnarray*} 
The formulas (\ref{es:s3}) can be generalized for interval $[a, b]$
other than [-1, 1] by the linear change of the variable \ 
$h(\tau)=\frac{1}{2}(b-a)\tau+\frac{1}{2}(a+b).$
Thus if \ 
$\eta_{j}=h(\tau_{j}), \hspace{3mm} j = 0,1,...,n,$ \ we have
\begin{eqnarray*}
\left[{\bf{I}}+\frac{b-a}{2}{\bf{K}}{\bf{D}}_{\sigma}\right]{\bar{\bf{x}}}={\bar{\bf{y}}},
\end{eqnarray*}
where
\begin{eqnarray*}
{\bf{K}} & = & (k(\eta_{i},\eta_{j}))^{n}_{i,j=0} , \ \ \
{\bf{D}}_{\sigma}  = diag(\sigma_{0}, \sigma_{1},...,\sigma_{n}),\\
{\bar{\bf{x}}} & = &
\left[x(\eta_{0}), x(\eta_{1}),...,x(\eta_{n})\right]^{T}, \ \
{\bar{\bf{y}}}  = 
\left[y(\eta_{0}), y(\eta_{1}),...,y(\eta_{n})\right]^{T}.\\
\end{eqnarray*} 
The accuracy of this discretization when $k$ and $x$ are not
 polynomials
 is discussed in a more general setting in
the next section.

\section{Gauss Type Quadrature for a Semismooth Kernel}

We consider now more general semismooth kernels, as in Definition 1,
for  which we write
\begin{equation}
\label{es:s4}
x(t)+\int^{t}_{a}k_{1}(t,s)x(s)ds+\int^{b}_{t}k_{2}(t,s)x(s)ds
=y(t), \hspace{5mm} a \leq t \leq b.
\end{equation}
In this section we describe the numerical technique for discretizing 
the equation (\ref{es:s4}).
It is based on the Clenshaw-Curtis quadrature described in Section 2, 
which is well suited for computing 
antiderivatives.
Let $$F(t)=\int^{t}_{-1}k_{1}(t,s)x(s)ds, \ \ \
 G(t,\lambda)=\int^{\lambda}_{-1}k_{1}(t,s)x(s)ds,$$
such that $F(t)=G(t,t),$ and let
 $$H(t)=\int^{1}_{t}k_{2}(t,s)x(s)ds, \ \ \
 J(t,\lambda)=\int^{1}_{\lambda}k_{2}(t,s)x(s)ds.$$
Further, assume that $k_{1}(t_{k},s)x(s)$ can be expanded in a finite set of polynomials,
i.e., $k_{1}(t_{k},s)x(s)=\sum^{n}_{i=0}\alpha_{ki}T_{i}(s).$
As we have seen in Section 2, if 
\begin{equation}
\label{es:s5}
G(t_{k},\lambda)=\sum_{j=0}^{n+1}\beta_{kj}T_{j}(\lambda),
\end{equation}
then 
\begin{equation}
\label{es:s6}
\left[\beta_{k0}, \ \beta_{k1},...,\beta_{kn+1}\right]^{T}={\bf{S}}_{L}
\left[\alpha_{k0},\alpha_{k1},...,\alpha_{kn}\right]^{T}.
\end{equation}
Similarly, assume that $k_{2}(t_{k},s)x(s)=\sum^{n}_{j=0}{\tilde{\alpha}}_{kj}T_{j}(s).$
If
\begin{eqnarray*}
J(t_{k},\lambda)  = \int^{1}_{\lambda}k_{2}(t_{k},s)x(s)ds
 =  \sum_{j=0}^{n+1}\tilde{\beta}_{kj}T_{j}(\lambda),
\end{eqnarray*}
then
\begin{equation}
\label{es:s7}
\left[\tilde{\beta}_{k0}, \ \tilde{\beta}_{k1},...,\tilde{\beta}_{kn+1}\right]^{T}
={\bf{S}}_{R}
\left[{\tilde{\alpha}}_{k0},{\tilde{\alpha}}_{k1},...,{\tilde{\alpha}}_{kn}\right]^{T},
\end{equation}
where 
\begin{eqnarray*}
{\bf{S}}_{R}=
\left[\begin{array}{llllll}
1 & 1 & 1 & 1 & \cdots & 1 \\
0 & -1 & 0  & 0 & \cdots & 0 \\
0 & 0 & -1 & 0 & \cdots & 0 \\
\vdots & \vdots & \vdots & \ddots & \ddots & \vdots \\
0 & 0 & 0 & 0 & -1& 0 \\
0 & 0 & \cdots & 0 & 0 & -1 \\
\end{array}
\right]
\left[\begin{array}{llllll}
0 & 0 & 0 & 0 & \cdots & 0 \\
1 & 0 & \frac{-1}{2}  & 0 & \cdots & 0 \\
0 & \frac{1}{4} & 0 & \frac{-1}{4} & \cdots & 0 \\
\vdots & \vdots & \ddots & \ddots & \ddots & \vdots \\
0 & \cdots & 0 & \frac{1}{2(n-1)} & 0 & \frac{-1}{2(n-1)} \\
0 & \cdots & 0 & 0 & \frac{1}{2n} & 0 \\
\end{array}
\right]
\end{eqnarray*}   
is the right spectral integration matrix.
Let $\tau_{k}, k=0,1,...,n,$ denote the zeros of $T_{n+1}.$
Substituting $\lambda=\tau_{k}, k=0,1,...,n,$ into (\ref{es:s5}), we obtain that
\begin{eqnarray*}
\left[\begin{array}{c}
G(t_{k},\tau_{0})\\
G(t_{k},\tau_{1})\\
\vdots\\
G(t_{k},\tau_{n})
\end{array}
\right]
={\bf{C}}{\bf{S}}_{L}{\bf{C}}^{-1}
diag(k_{1}(t_{k},\tau_{0}),...,k_{1}(t_{k},\tau_{n}))
\left[\begin{array}{c}
x(\tau_{0})\\
\vdots\\
x(\tau_{n})
\end{array}
\right]
\end{eqnarray*}
and, similarly,
\begin{eqnarray*}
\left[\begin{array}{c}
J(t_{k},\tau_{0})\\
J(t_{k},\tau_{1})\\
\vdots\\
J(t_{k},\tau_{n})
\end{array}
\right]
={\bf{C}}{\bf{S}}_{R}{\bf{C}}^{-1}
diag(k_{2}(t_{k},\tau_{0}),...,k_{2}(t_{k},\tau_{n}))
\left[\begin{array}{c}
x(\tau_{0})\\
\vdots\\
x(\tau_{n})
\end{array}
\right].
\end{eqnarray*}
We remark that in writing the equality sign in (\ref{es:s6}) and 
(\ref{es:s7}), we assume that
$\beta_{n+1}$ is set to zero. This is an acceptable assumption because
$G(t_{k},\lambda)$ is itself only approximately represented by the polynomials
in (\ref{es:s5}) and the overall accuracy is not affected.
Since $F(\tau_{k})=G(\tau_{k},\tau_{k})$ \  we get
\begin{eqnarray*}
F(\tau_{k})  & = &  \left[0,...,0,1,0,...,0\right]{\bf{C}}{\bf{S}}_{L}{\bf{C}}^{-1}
diag(k_{1}(\tau_{k},\tau_{0}),...,k_{1}(\tau_{k},\tau_{n}))
\left[\begin{array}{c}
x(\tau_{0})\\
\vdots\\
x(\tau_{n})
\end{array}
\right]\\
 & = &
\left[w_{k0},w_{k1},...,w_{kn}\right]
diag(k_{1}(\tau_{k},\tau_{0}),...,k_{1}(\tau_{k},\tau_{n}))
\left[\begin{array}{c}
x(\tau_{0})\\
\vdots\\
x(\tau_{n})
\end{array}
\right]\\
 & = & 
\left[w_{k0},w_{k1},...,w_{kn}\right]
diag(x(\tau_{0}),...,x(\tau_{n}))
\left[\begin{array}{c}
k_{1}(\tau_{k},\tau_{0})\\
\vdots\\
k_{1}(\tau_{k},\tau_{n})
\end{array}
\right]
\end{eqnarray*}
where $[w_{k0},...,w_{kn}]$ is the $(k+1)$-{st} row  of the matrix
${\bf{W}}\stackrel{\rm def}{=}{\bf{C}}{\bf{S}}_{L}{\bf{C}}^{-1}.$
We need now the following identity which can be verified by direct
calculation.
\begin{lemma}
Let ${\bf{A}}$ and \ ${\bf{B}}$ be $n \times n$ \ matrices and
 ${\bf{c}}=[c_{1},...,c_{n}]^{T}.$
Then $({\bf{A}} \circ
 {\bf{B}}) {\bf{c}}=diag({\bf{A}}diag(c_{1},...,c_{n})
{\bf{B}}^{T}),$
where ${\bf{A}} \circ {\bf{B}}$ denotes the Schur product of ${\bf{A}}$ and
$ {\bf{B}},$ \  $({\bf{A}} \circ {\bf{B}})_{ij}=a_{ij}b_{ij}, \ \
 i, j=1,...,n.$
\end{lemma}

Using this lemma we find that,
\begin{eqnarray}
\label{es:s8}
\left[\begin{array}{c}
F(\tau_{0})\\
F(\tau_{1})\\
\vdots\\
F(\tau_{n})
\end{array}
\right]
 =
diag({\bf{W}}diag(x(\tau_{0}),...,x(\tau_{n})){\bf{K}}_{1}^{T})
 =
({\bf{W}}\circ{\bf{K}}_{1})
\left[\begin{array}{c}
x(\tau_{0})\\
\vdots\\
x(\tau_{n})
\end{array}
\right],
\end{eqnarray}
where \ ${\bf{K}}_{1}=(k_{1}(\tau_{i},\tau_{j}))^{n}_{i,j=0}.$ 
Similarly,
\begin{eqnarray}
\label{es:s9}
\left[\begin{array}{c}
H(\tau_{0})\\
H(\tau_{1})\\
\vdots\\
H(\tau_{n})
\end{array}
\right]=({\bf{V}}\circ{\bf{K}}_{2})
\left[\begin{array}{c}
x(\tau_{0})\\
\vdots\\
x(\tau_{n})
\end{array}
\right],
\end{eqnarray}
where \ \ ${\bf{V}}={\bf{C}}{\bf{S}}_{R}{\bf{C}}^{-1}.$
The formulas (\ref{es:s8}) and (\ref{es:s9}) can be generalized 
for an interval $[a, b]$
other than $[-1, 1]$ by the linear change of variables,
 $h(\tau)=\frac{1}{2}(b-a)\tau+\frac{1}{2}(a+b).$
Thus if
$\eta_{j}=h(\tau_{j}), \hspace{3mm} j = 0,1,...,n,$
and with the notation
\begin{eqnarray*}
F_{a}(t)=\int^{t}_{a}k_{1}(t,s)x(s)ds, \ \ \
H_{b}(t)=\int^{b}_{t}k_{2}(t,s)x(s)ds,
\end{eqnarray*}
we have
\begin{eqnarray}
\label{es:s10}
\left[\begin{array}{c}
F_{a}(\eta_{0})\\
F_{a}(\eta_{1})\\
\vdots\\
F_{a}(\eta_{n})
\end{array}
\right]=
\frac{b-a}{2}
({\bf{W}}\circ{\bf{K}}_{1})
\left[\begin{array}{c}
x(\eta_{0})\\
x(\eta_{1})\\
\vdots\\
x(\eta_{n})
\end{array}
\right]
\end{eqnarray}
and, 
\begin{eqnarray}
\label{es:s11}
\left[\begin{array}{c}
H_{b}(\eta_{0})\\
H_{b}(\eta_{1})\\
\vdots\\
H_{b}(\eta_{n})
\end{array}
\right]=\frac{b-a}{2}({\bf{V}}\circ{\bf{K}}_{2})
\left[\begin{array}{c}
x(\eta_{0})\\
x(\eta_{1})\\
\vdots\\
x(\eta_{n})
\end{array}
\right].
\end{eqnarray}
Using (\ref{es:s10}) and (\ref{es:s11}) we can now discretize 
the equation (\ref{es:s4}) as follows,
\begin{equation}
\label{es:s12}
\left[
{\bf{I}}+\frac{b-a}{2}({\bf{W}}\circ{\bf{K}}_{1}+{\bf{V}}\circ{\bf{K}}_{2})
\right]{\bar{\bf{x}}}={\bar{\bf{y}}},
\end{equation}
where \hspace{3mm} ${\bar{\bf{x}}}=[x(\eta_{0}),...,x(\eta_{n})]^{T}$
and \ ${\bar{\bf{y}}}=[y(\eta_{0}),...,y(\eta_{n})]^{T}.$
Next we show that if \
the kernel function $k(t,s)$ is smooth, such that $k_1=k_2$,
  then the discretization
 (\ref{es:s12}) reduces to (\ref{es:s3}). 
\begin{prop}
Suppose that $k(t,s) \in C^{p}_{[a,b]\times [a,b]}$, 
and that $k_{1}(t,s)=k_{2}(t,s)=k(t,s)$. 
Then, 
\begin{eqnarray*}
\left[{\bf{I}}+\frac{b-a}{2}({\bf{W}}\circ{\bf{K}}_{1}+{\bf{V}}\circ{\bf{K}}_{2})\right]
{\bar{\bf{x}}}=
\left[{\bf{I}}+\frac{b-a}{2}{\bf{K}}{\bf{D}}_{\sigma}\right]{\bar{\bf{x}}}.
\end{eqnarray*}.
\end{prop}
\bpr
Without any loss of generality we assume that $a=-1$ \ and $b=1.$
For $t = t_{k}$, \hspace{3mm}
$-1 \leq t_{k} \leq 1,$ 
\begin{equation}
\label{es:s13}
x(t_{k})+G(t_{k},\lambda)+J(t_{k},\lambda)=y(t_{k})
\end{equation}
holds for  any $\lambda$, \ $-1 \leq \lambda \leq 1.$ \ \
Therefore if $\lambda=1$, then  $J(t_{k},1)=0$ and
\begin{eqnarray*}
G(t_{k},1) & = & \sum^{n}_{j=0}\beta_{kj} \
 =  [1,...,1]{\bf{S}}_{L}
\left[\begin{array}{c}
\alpha_{k0}\\
\vdots\\
\alpha_{kn}
\end{array}
\right] \\
& = &
[1,...,1]{\bf{S}}_{L}
{\bf{C}}^{-1}
diag(k(t_{k},\tau_{0}),...,k(t_{k},\tau_{n}))
\left[\begin{array}{c}
x(\tau_{0})\\
\vdots\\
x(\tau_{n})
\end{array}
\right] \\
& = &
[\sigma_{0}, \sigma_{1},...,\sigma_{n}]
diag(k(t_{k},\tau_{0}),...,k(t_{k},\tau_{n}))
\left[\begin{array}{c}
x(\tau_{0})\\
\vdots\\
x(\tau_{n})
\end{array}
\right]\\
& = &
[k(t_{k},\tau_{0}),k(t_{k},\tau_{1}),...,k(t_{k},\tau_{n})]
diag(\sigma_{0}, \sigma_{1},...,\sigma_{n})
\left[\begin{array}{c}
x(\tau_{0})\\
x(\tau_{1})\\
\vdots\\
x(\tau_{n})
\end{array}
\right].
\end{eqnarray*}
Substituting $t_{k}=\tau_{k}$ \ for $k=0,1,...,n,$ into (\ref{es:s13}),
 we obtain that
$$\left[{\bf{I}}+{\bf{K}}{\bf{D}}_{\sigma}\right]{\bar{\bf{x}}}={\bar{\bf{y}}}.$$
The assertion now follows.
\epr

We compared the numerical behavior of the discretization (\ref{es:s3}) and 
(\ref{es:s12})
for a number of smooth kernels and found that numerical answers
differed in  accuracy at the level of  machine precision only.

We now estimate the accuracy of approximation of the integral equation
(\ref{es:s4}) with the linear system of equations (\ref{es:s12}). 
The following property
of Chebyshev expansions can be derived along the lines
of an argument in  Gottlieb and Orszag
(\cite{Got}, p.29). 
\begin{prop}
Let $f\in C^{r}[-1, 1], \ \ r > 1,$ \ \ and let
$$f(t)=\sum^{\infty}_{j=0}\alpha_{j}T_{j}(t), \hspace{3mm} -1 \leq t \leq 1.$$
Then
$$|\alpha_{j}| \leq \frac{2}{\pi}\int^{\pi}_{0}\left|\frac{d^{r}}{d\theta^{r}}
f(cos\theta)\right|d\theta\frac{1}{j^{r}}=\frac{c}{j^{r}}$$
and
$$\left|f(t)-\sum^{n}_{j=0}\alpha_{j}T_{j}(t)\right|\leq\frac{c}{r-1}
\frac{1}{n^{r-1}}.$$
\end{prop}
It implies that if $f(r)$ is
analytic then the convergence of  Chebyshev expansions 
is superalgebraic. Let now
$F_{l}(x)=\int^{x}_{-1}f(t)dt$ \hspace{3mm} and \hspace{3mm}
$F_{r}(x)=\int^{1}_{x}f(t)dt$.
The following result can be found in Greengard and Rokhlin \cite{Green}.
\begin{prop} 
Suppose that $f \in C^{r}_{[-1,1]}, \ \ r>1,$ 
and that ${\bar{f}}=(f(\tau_{0}),...,f(\tau_{n}))^{T},$  
is the vector of  the function values at the roots
of \  $T_{n+1}(x)$. 
Suppose further that ${\bar{F}}_{l}$ \ \ and ${\bar{F}}_{r}$ are defined by \\
$${\bar{F}}_{l}=(F_{l}(\tau_{0}),...,F_{l}(\tau_{n}))^{T}, \ \ \  
{\bar{F}}_{r}=(F_{r}(\tau_{0}),...,F_{r}(\tau_{n}))^{T}.$$
Then
$$||{\bar{F}}_{l}-{\bf{C}}{\bf{S}}_{L}{\bf{C}}^{-1}{\bar{f}}||_{\infty}
={\it{O}}(\frac{1}{n^{r-1}})$$
and
$$||{\bar{F}}_{r}-{\bf{C}}{\bf{S}}_{R}{\bf{C}}^{-1}{\bar{f}}||_{\infty}=
{\it{O}}(\frac{1}{n^{r-1}}).$$
Furthermore, all elements of the matrix ${\bf{C}}{\bf{S}}_{L}{\bf{C}}^{-1}$
and $ {\bf{C}}{\bf{S}}_{R}{\bf{C}}^{-1} $ are strictly positive.
\end{prop}

Let now \
${\eta}_{i}=\frac{b-a}{2}{\tau}_{i}+\frac{a+b}{2},$
where \ $\tau_{i}$ is a zero of \  $T_{n+1}(x),$
 \ \ for \ \ $i=0,1,...,n,$ be the shifted Chebyshev points,
and \ ${\hat{\bf{x}}}=(x({\eta}_{0}),x({\eta}_{1}),...,x({\eta}_{n}))^{T}$
\ be the vector of values of solution $x(t)$ of equation (\ref{es:s4})
 at ${\eta}_{i}.$
The following proposition follows immediately from standard properties 
of the Riemann integral, (see e.g. \cite{Royden}, p.105).  
\begin{prop}
Let $k(t,s)$ be $(p_{1},p_{2})$-semismooth and
let  $y(t) \in C^{q}_{[a,b]},$ \ such that $r=\min \{p_{1},p_{2},q\}>1.$
Let the equation (\ref{es:s4}) define an invertible operator on $C^{r}_{[a,b]}.$
Then $x \in C^{r}_{[a,b]}.$
\end{prop}

Let now
$${\bar{F}}_{a}=(F_{a}(\eta_{0}),...,F_{a}(\eta_{n}))^{T}$$
and $${\bar{H}}_{b}=(H_{b}(\eta_{0}),...,H_{b}(\eta_{n}))^{T}.$$
It follows from Proposition 3 that
in conditions of Proposition 4,
$$||{\bar{F}}_{a}-\frac{b-a}{2}({\bf{W}}\circ {\bf{K}}_{1})
{\hat{\bf {x}}}||_{\infty}
={\it{O}}(\frac{1}{n^{r-1}}),$$
and
$$||{\bar{H}}_{b}-\frac{b-a}{2}({\bf{V}}\circ {\bf{K}}_{2})
{\hat{\bf {x}}}||_{\infty}
={\it{O}}(\frac{1}{n^{r-1}}).$$
Combining the above results we obtain the following estimate
for the residual.
\begin{theorem}
Let ${\bar {\bf{x}}}$ be a solution vector of the
equation (\ref{es:s12}), and ${\hat {\bf {x}}}$ the vector of values
of the solution $x(t)$ at $t={\eta}_{i}, \ i=0,1,...,n.$
Suppose that $k(t,s)$ is $(p_{1},p_{2})$-semismooth, and that 
$y(t) \in C^{q}_{[a,b]}.$ 
Suppose further that the equation (\ref{es:s4}) defines an invertible
operator on $C^{r}_{[a,b]},$  where $r=\min \{p_{1},p_{2},q\}>1.$
Then,
$$||({\bf{I}}+\frac{b-a}{2}({\bf{W}}\circ {\bf{K}}_{1}+{\bf{V}}
\circ {\bf{K}}_{2}))({\hat {\bf{x}}}-{\bar {\bf{x}}})||_{\infty} 
={\it{O}}(\frac{1}{n^{r-1}}).$$ 
\end{theorem}

It follows from the collectively compact operator theory, see  Anselone 
\cite{Ans}, that for sufficiently large $n$
the matrices ${\bf{I}}+\frac{b-a}{2}({\bf{W}}\circ {\bf{K}}_{1}+
{\bf{V}}\circ {\bf{K}}_{2})$, which depend on $n$,
 are invertible and their inverses are
uniformly bounded. Therefore Theorem 1 implies that for increasing
$n,$  the convergence of  ${\bar{x}}$  to ${\hat{x}}$ is 
of order ${\it{O}}({n^{1-r}}).$
If $p_{1}=p_{2}=q=\infty,$  then the convergence is
superalgebraic. Numerical examples in Section 5 indeed demonstrate
this type of
convergence.

\section{The Composite Rule}
In this section, we describe the composite rule corresponding to
the quadrature of (\ref{es:s12}).                                              Let $$a=b_{0} \leq b_{1} \leq \cdots \leq b_{m}=b$$
be a partition of the interval $[a,b]$,  and let
$${\tau}^{(j)}_{k}=\frac{1}{2}(b_{j}-b_{j-1}){\tau_{k}}+
\frac{1}{2}(b_{j}+b_{j-1}), \ \ \ k=0,1,...,n_j,$$
be the Chebyshev support points mapped into $[b_{j-1},b_{j}].$ 
Define 
\[
x(t)=
\left\{ \begin{array}{cc}
x_{1}(t) & {\rm if} \ \      b_{0} \leq t \leq b_{1} \\
x_{2}(t) & {\rm if} \ \ \ \  b_{1} < t \leq b_{2} \\
\vdots \\
x_{m}(t) & {\rm if} \ \ \ \  b_{m-1} < t \leq b_{m},
\end{array}
\right.
\]
and
\[
y(t)=
\left\{ \begin{array}{cc}
y_{1}(t) & {\rm if} \ \     b_{0} \leq t \leq b_{1}\\
y_{2}(t) & {\rm if} \ \ \ \ b_{1} < t \leq b_{2} \\
\vdots \\
y_{m}(t) & {\rm if} \ \ \ \ b_{m-1} < t \leq b_{m},
\end{array}
\right.
\] 
and rewrite the equation (\ref{es:s1}) as a system of $m$
equations, for $j=1,...,m$,
\begin{eqnarray*}
x_{j}(t)  +  \int^{b_{1}}_{b_{0}}k_{1}(t,s)x_{1}(s)ds+ \cdots +
\int^{t}_{b_{j-1}}k_{1}(t,s)x_{j}(s)ds 
  +  \int^{b_{j}}_{t}k_{2}(t,s)x_{j}(s)ds 
\end{eqnarray*}
\begin{equation}
\label{es:s14} 
+ \cdots +
\int^{b_{m-1}}_{b_{m}}k_{2}(t,s)x_{m}(s)ds
=y_{j}(t). 
\end{equation}
Applying the quadrature of (\ref{es:s12}) to each of the integrals
we obtain a system of linear equations as follows, for $j=1,...,m$,
\begin{eqnarray*}
\frac{b_{1}-b_{0}}{2}[({\bf{W}}+{\bf{V}}) \circ
{\bf{K}}_{1j}]{\bar{x}}_{1}+ \cdots   
+[{\bf{I}}+\frac{b_{j}-b_{j-1}}{2}({\bf{W}}\circ {\bf{K}}_{jj}+
{\bf{V}}\circ {\tilde{\bf{K}}}_{jj})]
{\bar{x}}_{j}+  \\ 
 + \cdots +\frac{b_{m}-b_{m-1}}{2}[({\bf{W}}+{\bf{V}}) 
\circ {\bf{K}}_{mj}]{\bar{x}}_{m}
={\bar{y}}_{j},
\end{eqnarray*}
where
\begin{eqnarray*}
{\bar{x}}_{j} & = & [x({\tau}^{(j)}_{0}),x({\tau}^{(j)}_{1}),...,
x({\tau}^{(j)}_{n_j})]^{T},
\hspace{3mm}
{\bar{y}}_{j}=[y({\tau}^{(j)}_{0}),y({\tau}^{(j)}_{1}),...,
y({\tau}^{(j)}_{n_j})]^{T}, \\ 
{\bf{K}}_{jj} & = & (k_{1}({\tau}^{(j)}_{p},{\tau}^{(j)}_{q}))_{p,q=0}^{n_j}, \hspace{5mm} 
{\tilde{\bf{K}}}_{jj}=(k_{2}({\tau}^{(j)}_{p},{\tau}^{(j)}_{q}))_{p,q=0}^{n_j}, \\ 
{\bf{K}}_{ij} & = & (k_{1}({\tau}^{(j)}_{p},{\tau}^{(i)}_{q}))_{p,q=0}^{n_j,n_i},\ \ \ \ 
if \ \ \  i<j, \\
{\bf{K}}_{ij} & = & (k_{2}({\tau}^{(j)}_{p},{\tau}^{(i)}_{q}))_{p,q=0}^{n_j,n_i},\ \ \ \ 
if \ \ \  i>j, \\
\end{eqnarray*}
or in a block matrix form,
\begin{eqnarray}
\label{es:s15}
\left[\begin{array}{cccc}
{\bf{A}}_{11} &  {\bf{A}}_{12} & \cdots & {\bf{A}}_{1m} \\
{\bf{A}}_{21} &  {\bf{A}}_{22} & \cdots & {\bf{A}}_{2m} \\
\vdots \\
{\bf{A}}_{m1} &  {\bf{A}}_{m2} & \cdots & {\bf{A}}_{mm}
\end{array}
\right]
\left[\begin{array}{c}
{\bar{x}}_{1} \\
{\bar{x}}_{2} \\
\vdots \\
{\bar{x}}_{m}
\end{array}
\right]
=
\left[\begin{array}{c}
{\bar{y}}_{1} \\
{\bar{y}}_{2} \\
\vdots \\
{\bar{y}}_{m}
\end{array}
\right],
\end{eqnarray}
where
\begin{eqnarray*}
{\bf{A}}_{jj} & = & [{\bf{I}}+\frac{b_{j}-b_{j-1}}{2}({\bf{W}}\circ 
{\bf{K}}_{jj}+{\bf{V}}\circ {\tilde{\bf{K}}}_{jj})],\\
{\bf{A}}_{ij} & = & \frac{b_{j}-b_{j-1}}{2}[({\bf{W}}+{\bf{V}}) 
\circ {\bf{K}}_{ji}], \ \ \ \  if \ \ \ i \neq j. 
\end{eqnarray*}

In this paper we do not consider the issue of how to partition the interval 
$[a,b]$. An adaptive quadrature rule 
is possible here along the same lines as in 
\cite{reo2}, \cite{GR}, 
namely, by using the size of last Chebyshev coefficients
of $k_1, k_2$ and $y$ in a given subinterval
of partition to determine whether this subinterval
 should be further subdivided.
This adaptive rule is a part of our research project
in which we are going to compare the algorithm of the present paper 
with existing algorithms for Schroedinger equations with 
physically realistic nonlocal potentials. 

In general, the matrix (\ref{es:s15}) is not structured,
and is being solved by the standard Gaussian elimination at the cost of 
$O(m^3)$ arithmetic operations, (we assume here that $m$ is much larger
than the $n_j$'s). If however, the semismooth 
kernel $k(t,s)$ has some additional structure, then this 
structure is usually inherited by the matrix in (\ref{es:s15}). 
For example, if 
$k_1$ and $k_2$ are low rank kernels then the matrix $A$ becomes semiseparable,
and can be solved by existing linear complexity algorithms.
We remark that in the case of the Schroedinger
equation with non-local potentials
discussed in Section 7 below,
the overall kernel is obtained as the composition
of the semi-separable Green's function with the non-local
potential. 
If the non-local potential is also semi-separable, 
which is the case when
the non-locality arises from exchange terms, then the overall
kernel is semi-separable as well, and the numerical techniques
presented here, although still applicable, can be replaced by
the IEM methods (\cite {reo}, \cite {reo2}) previously developed
for local potentials. These methods give highly accurate linear complexity 
algorithms for the integral equation itself.

If the kernel $k(t,s)$ depends on the difference of the arguments,
\[
k(t,s)=k(|t-s|)=
\left \{ \begin{array}{cc}
k_{1}(t-s) & {\rm if} \ \  0 \leq s \leq t  \\
k_{2}(t-s) & {\rm if} \ \ \ \  t < s \leq  T , 
\end{array}
\right.
\]
and if we use a uniform partition with the same number of points
per partition, then, 
\begin{eqnarray*}
k_{r}({\tau}^{(i)}_{p},{\tau}^{(i)}_{q}) = k_{r}({\tau}_{p}-{\tau}_{q}),
 \ \ \   
r=1,2,
\end{eqnarray*}
and we obtain a block Toeplitz matrix,
\begin{eqnarray}
\label{es:s16}
\left[ \begin{array} {ccccc}
{\bf{A}}_{1} & {\tilde{\bf{A}}}_{2} & {\tilde{\bf{A}}}_{3} & \cdots
& {\tilde{\bf{A}}}_{m} \\
{\bf{A}}_{2} & {\bf{A}}_{1} & {\tilde{\bf{A}}}_{2} & \cdots &
{\tilde{\bf{A}}}_{m-1} \\  
{\bf{A}}_{3} & {\bf{A}}_{2} & {\bf{A}}_{1} & \cdots &
{\tilde{\bf{A}}}_{m-2} \\
       & \ddots & \ddots & \ddots &        \\
{\bf{A}}_{m} & {\bf{A}}_{m-1} & \cdots & {\bf{A}}_{2} & {\bf{A}}_{1}
\end{array}
\right]
\left[ \begin{array}{c}
{\bar{x}}_{1} \\
{\bar{x}}_{2} \\
\vdots \\
{\bar{x}}_{m}
\end{array}
\right]
=\left[ \begin{array}{c}
{\bar{y}}_{1} \\
{\bar{y}}_{2} \\
\vdots \\
{\bar{y}}_{m}
\end{array}
\right],
\end{eqnarray}
which can be solved in $O(m \log (m))$ or $O(m^2)$ arithmetic operations by 
algorithms available in the literature.

Finally we would like to point out that the above composite rule can be used 
to handle kernels which have a finite number of singularities
on the main diagonal, $t=s$.
In this case one has to 
include all the singular points as a subset of all partition points.
This is illustrated in the next section on a simplest
case of one such singularity, but in full detail. 

\section{Kernels with Singularities on the Main Diagonal}

Suppose that the kernel $k(t,s)$ has a singularity at $ (c, c), 
a < c < b,$ \  inside the square 
$[a, b] \times [a, b]$.  \  Since the Chebyshev points 
${\eta}_{i}=\frac{b-a}{2}{\tau}_{i}+\frac{a+b}{2}, \  
i=0,1,2,...,n,$ are clustered towards the boundaries of the
square $[a,b]\times [a,b]$, we do not have 
sufficient values of $k(t,s)$ around  singular point $(c,c)$
 for a good approximation.
Therefore we apply the composite rule of Section 4 
with $c$ being a partition point.
For the sake of simplicity we assume that $c$ is the only partition point,
thus $[a,b]$ is partitioned 
 into two subintervals $[a,c]$ and $[c,b].$
Without loss of generality, we consider the solution of
\begin{equation}
\label{es:s17}
x(t)+\int^{1}_{-1}k(t,s)x(s)ds=y(t),
\end{equation} 
where
\[
k(t,s)=
\left \{ \begin{array}{cc}
k_{1}(t,s) & {\rm if} \ \  -1 \leq s \leq t  \\
k_{2}(t,s) & {\rm if} \ \ \ \  t < s \leq 1 , 
\end{array}
\right.
\] 
and assume that the kernel $k(t,s)$ has a singular point at (0,0).
Define
\[
x(t)=
\left\{ \begin{array}{cc}
x_{1}(t) & {\rm if} \ \  -1 \leq t \leq 0 \\
x_{2}(t) & {\rm if} \ \ \ \  0 < t \leq 1
\end{array}
\right.
\]
and
\[
y(t)=
\left\{ \begin{array}{cc}
y_{1}(t) & {\rm if} \ \ -1 \leq t \leq 0 \\
y_{2}(t) & {\rm if} \ \ \ \  0 < t \leq 1 .
\end{array}
\right.
\] 
We can rewrite the equation (\ref{es:s17}) as a system of two equations. 
For $ -1 \leq t \leq 0,$ 
\begin{equation}
\label{es:s18}
x_{1}(t)+\int^{t}_{-1}k(t,s)x_{1}(s)ds+\int^{0}_{t}k(t,s)x_{1}(s)ds
+\int^{1}_{0}k(t,s)x_{2}(s)ds=y_{1}(t),
\end{equation}
and for $0 < t \leq 1,$
\begin{equation}
\label{es:s19}
x_{2}(t)+\int^{0}_{-1}k(t,s)x_{1}(s)ds+\int^{t}_{0}k(t,s)x_{2}(s)ds
+\int^{1}_{t}k(t,s)x_{2}(s)ds=y_{2}(t).
\end{equation}
Discretizations of  the equations (\ref{es:s18}) and (\ref{es:s19}), 
respectively,
are as follows.
\begin{eqnarray}
\label{es:s20}
[{\bf{I}}+\frac{1}{2}({\bf{W}}\circ {\bf{K}}_{11}+{\bf{V}}\circ 
{\tilde{\bf{K}}}_{11})]
{\bar{x}}_{1}+\frac{1}{2}[({\bf{W}}+{\bf{V}}) \circ {\bf{K}}_{21}]{\bar{x}}_{2}
={\bar{y}}_{1}
\end{eqnarray}   
and
\begin{eqnarray}
\label{es:s21}
\frac{1}{2}[({\bf{W}}+{\bf{V}}) \circ {\bf{K}}_{12}]{\bar{x}}_{1}+
[{\bf{I}}+\frac{1}{2}({\bf{W}}\circ {\bf{K}}_{22}+{\bf{V}}\circ 
{\tilde{\bf{K}}}_{22})]
{\bar{x}}_{2}={\bar{y}}_{2},
\end{eqnarray}
where
\begin{eqnarray*}
{\bf{K}}_{11}=k_{1}(\tau_{p}^{(1)},\tau_{q}^{(1)})_{p,q=0}^{n_1}, \hspace{5mm} 
{\tilde{\bf{K}}}_{11}=k_{2}(\tau_{p}^{(1)},\tau_{q}^{(1)})_{p,q=0}^{n_1}, 
\hspace{5mm} {\bf{K}}_{12}=k_{1}(\tau_{p}^{(2)},\tau_{q}^{(1)})
_{p,q=0}^{n_2,n_1}, \\
{\bf{K}}_{22}=k_{1}(\tau_{p}^{(2)},\tau_{q}^{(2)})
_{p,q=0}^{n_2}, \hspace{5mm} 
{\tilde{\bf{K}}}_{22}=k_{2}(\tau_{p}^{(2)},\tau_{q}^{(2)})_{p,q=0}^{n_2},
\hspace{5mm} {\bf{K}}_{21}=k_{2}(\tau_{p}^{(1)},\tau_{q}^{(2)})
_{p,q=0}^{n_1,n_2},
\end{eqnarray*}
with
\begin{eqnarray*}
\tau_{i}^{(1)}=\frac{1}{2}\tau_{i}-\frac{1}{2} \hspace{4mm} and 
\hspace{4mm} \tau_{j}^{(2)}=\frac{1}{2}\tau_{j}+\frac{1}{2},
\hspace{3mm} i=0,1,...,n_{1}, \hspace{3mm} j=0,1,...,n_{2}.
\end{eqnarray*}
Here ${\bar{x}}_{1}=[x_{1}(\tau_1^{(1)}),...,x_{1}(\tau_{n_1}^{(1)})]^{T}$ 
\hspace{3mm} and \hspace{3mm} 
${\bar{x}}_{2}=[x_{2}(\tau_1^{(2)}),...,x_{2}(\tau_{n_2}^{(2)})]^{T}.$ 
In the matrix form we write,
\begin{eqnarray}
\label{es:s22}
\left[\begin{array}{c}
{\bf{A}}_{11} \ \  {\bf{A}}_{12} \\
{\bf{A}}_{21} \ \  {\bf{A}}_{22}
\end{array}
\right]
\left[\begin{array}{c}
{\bar{x}}_{1} \\
{\bar{x}}_{2}
\end{array}
\right]
=
\left[\begin{array}{c}
{\bar{y}}_{1} \\
{\bar{y}}_{2}
\end{array}
\right],
\end{eqnarray}
where
\begin{eqnarray*}
{\bf{A}}_{11} & = & [{\bf{I}}+\frac{1}{2}({\bf{W}}\circ {\bf{K}}_{11}
+{\bf{V}}\circ {\tilde{\bf{K}}}_{11})],\\
{\bf{A}}_{12} & = & \frac{1}{2}[({\bf{W}}+{\bf{V}}) \circ {\bf{K}}_{21}],\\
{\bf{A}}_{21} & = & \frac{1}{2}[({\bf{W}}+{\bf{V}}) \circ {\bf{K}}_{12}],\\
{\bf{A}}_{22} & = & [{\bf{I}}+\frac{1}{2}({\bf{W}}\circ {\bf{K}}_{22}
+{\bf{V}}\circ {\tilde{\bf{K}}}_{22})].
\end{eqnarray*}
The size of the matrix in the equation (\ref{es:s22}) is  $(n_1+n_2+2) \times
(n_1+n_2+2)$. 
The formulas (\ref{es:s18}) and (\ref{es:s19}) can be generalized for 
interval $[a,c]$ and
$[c,b]$ other than
$[-1,1]$ by the linear change of the variable
$ p(t)=\frac{1}{2}(c-a)t+\frac{1}{2}(c+a)$ \ and \
$q(t)=\frac{1}{2}(b-c)t+\frac{1}{2}(b+c)$ 
if $(c, c)$  is a singular point
of $k(t,s)$.
Thus if 
$${\bar{\tau}}_{i}^{(1)}=p(\tau_{i}) \  and  \ \ \
  {\bar{\tau}}_{i}^{(2)}=q(\tau_{i}),  \hspace{3mm} i=0,1,2,...,n_i,$$
then,
\begin{eqnarray*}
\left[\begin{array}{c}
{\bar{\bf{A}}}_{11} \ \  {\bar{\bf{A}}}_{12} \\
{\bar{\bf{A}}}_{21} \ \  {\bar{\bf{A}}}_{22}
\end{array}
\right]
\left[\begin{array}{c}
{\bar{x}}_{1} \\
{\bar{x}}_{2}
\end{array}
\right]
=
\left[\begin{array}{c}
{\bar{y}}_{1} \\
{\bar{y}}_{2}
\end{array}
\right],
\end{eqnarray*}
where
\begin{eqnarray*}
{\bar{\bf{A}}}_{11} & = & [{\bf{I}}+\frac{c-a}{2}({\bf{W}}\circ {\bf{K}}_{11}
+{\bf{V}}\circ {\tilde{\bf{K}}}_{11})],\\
{\bar{\bf{A}}}_{12} & = & \frac{b-c}{2}[({\bf{W}}+{\bf{V}}) \circ {\bf{K}}_{21}],\\
{\bar{\bf{A}}}_{21} & = & \frac{c-a}{2}[({\bf{W}}+{\bf{V}}) \circ {\bf{K}}_{12}],\\
{\bar{\bf{A}}}_{22} & = & [{\bf{I}}+\frac{b-c}{2}({\bf{W}}\circ {\bf{K}}_{22}
+{\bf{V}}\circ {\tilde{\bf{K}}}_{22})],
\end{eqnarray*}
with $K_{i,j}$ defined as above.
A numerical example illustrating this technique is given in the next section.

\section{Numerical Examples}
In this section we compare our methods with
  some exiting algorithms for the following types of kernels, 
\begin{itemize}
\item{Type 1} :  Discontinuity  along the diagonal $t = s.$ 
\item{Type 2} :  Discontinuity in the first order partial derivatives along  
                 the diagonal $t = s.$ 
\item{Type 3} :  Singularity on the boundary of the square and Type 2. 
\item{Type 4} :  Singularity on the main diagonal.
\end{itemize}
These are the methods which have been implemented for comparison purposes,
\begin{description}
\item[G-Leg] : $Nystr{\ddot{o}}m$ discretization based on 
               Gauss-Legendre quadrature.
\item[T-Def] :  Two step Deferred Approach to the limit.
              Approximate solutions $x_{1},$ \ $x_{2}$,
             \ and \ $x_{3}$ 
             for subintervals of partition  $h,$ \ $\frac{h}{2}$,
             \ and \ $\frac{h}{4}$, respectively, are computed. Then the 
             numerical solution $x(s)$ is obtained by
              ( see e.g. Baker \cite{Baker})
             $$x(s)=\frac{64x_{3}(s)+x_{1}(s)-20x_{2}(s)}{45}.$$ 
           
\item[Atk-T] : Atkinson's iteration with the composite Trapezium rule. 
              Applied to  kernels $k(t,s)$ which have 
              discontinuities in the first
             order partial derivatives along the diagonal  $ t=s $.
\item[Alg-1] : Algorithm of Section 2, (\ref{es:s3}).
\item[Schur] : Algorithm of Section 3, (\ref{es:s12}).
\item[Sch-C] : Algorithm of Section 4, (\ref{es:s15}).
\end{description}

The number of points used in discretizations is denoted by $n.$
{\it Error} denotes $||x-x_{\tau}||/||x||$, where $x$ and $x_{\tau}$ are
the analytic and the numerical soluions, respectively.
In each plot, $log(Error)$ is the  common logarithm of the {\it Error}.
All computations were done on a DELL Workstation 
with operating
system RedHat Linux 5.2 in double precision. All examples are set-up by
choosing a simple analytic solution and then computing the
corresponding right hand side.
We remark that the values of $x(t)$ are found inside
the interval (or each of the subintervals of partition ) at
Chebyshev points ${\tau}_{0}, {\tau}_{1},..., {\tau}_{n}.$
The value of $x(t)$ for $t\neq {\tau}_{k}$ can be found
as follows. Applying ${\bf{C}}^{-1}$ we can find 
``Chebyshev-Fourier'' coeffcients of $x(t),$ 
\begin{eqnarray*}
\left[\begin{array}{c}
\alpha_{0}\\
\alpha_{1}\\
\vdots\\
\alpha_{n}
\end{array}
\right]={\bf{C}}^{-1}
\left[\begin{array}{c}
x({\tau}_{0})\\
x({\tau}_{1})\\
\vdots\\
x({\tau}_{n})
\end{array}
\right].
\end{eqnarray*}
Thus,
$$x(t) \cong \sum^{n}_{j=0}{\alpha}_{j}T_{j}(h(t)), \hspace{5mm} a \leq t \leq b.$$
The value of $T_{j}(t)$ for $t\neq {\tau}_{k}$ is found now
using the recursion satisfied by Chebyshev polynomials,
$T_{j+1}(t)=2tT_{j}(t)-T_{j-1}(t).$ \\
{\bf{Example}} 1. 
\begin{eqnarray*}
x(t)+\lambda \int^{1}_{-1}k(t,s)x(s)ds=y(t), \ \ \ -1 \leq t \leq 1,
\end{eqnarray*} 
where $y(t)=\lambda(e+e^{-1})+(1-2\lambda)e^{-t},$ \
and \\
\[
k(t,s)=
\left\{ \begin{array}{cc}
 1 & {\rm if} \ \ -1 \leq s \leq t \\
-1 & {\rm if} \ \ \ \  t < s \leq 1.
\end{array}
\right.
\]
The analytical solution is $x(t)=e^{-t}$. 
Since this kernel is discontinuous along the diagonal $ t = s $, 
Gauss-Legendre quadrature gives low accuracy. The accuracy in the
Atkinson's iteration improves very slowly.
The algorithm of Section 3 gives accuracy of order $10^{-15}$  with only 
16 support points, whereas the 2-step method of
Deferred Approach to the Limit method requires $n= 256$ points to
achieve comparable accuracy. Moreover, it requires computation of
$x_{2}(t)$ and $x_{3}(t)$ at the cost of ${\it{O}}((2n)^{3})$ and
${\it{O}}((4n)^{3}),$  respectively.
\begin{figure}[ht]
\begin{center}
\vspace{0in}
\leavevmode
\hspace{0in}
\epsfxsize =3in
\epsfbox{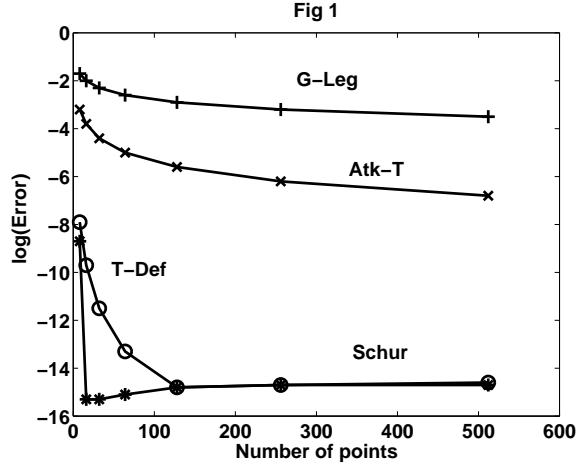}
\end{center}
\vspace{0in}
\hspace{0in}
\caption{ \label{figure1} Comparison of numerical solutions of Example 1,  \ \ $\lambda=0.1$. \ \   }
\end{figure}
\vspace{0.2in}
\\
{\bf{Example}} 2.
\begin{eqnarray*}
x(t)+{\lambda} \int^{T}_{0}k(t,s)x(s)ds=y(t),
 \ \ \ 0 \leq t \leq T
\end{eqnarray*} 
where   $y(t)=(1-\frac{{\lambda}sin^{2}(T)}{2}+{\lambda})sint+
(\frac{T}{2}-t-\frac{sin(2T)}{4}){\lambda}cost,$ \
and \\
\[
k(t,s) =  sin(|t-s|) \\
       = 
\left\{ \begin{array}{cc}
sin(t-s) & {\rm if} \ \  0 \leq s \leq t\\
sin(s-t) & {\rm if} \ \  t < s \leq \frac{\pi}{2}.
\end{array}
\right.
\]
The analytical solution is $x(t)=sin(t)$.
This  kernel has discontinuities in the first order partial
derivatives along the diagonal, $ t = s.$ 
Again standard $Nystr{\ddot{o}}m$-type
discretization methods  fail to give high accuracy in this case. 
In the first experiment
 we take $T=\frac{\pi}{2}$ \ and $\lambda=-\frac{4}{\pi}.$
Our method shows the order of $10^{-14}$ accuracy with only 16
points in $[0,\frac{\pi}{2}]$ without any partitioning.
The 2-step method of  Deferred Approach to the Limit  gives 
the accuracy of ${\it {O}}(10^{-14})$ with $n=256,$ 
but  at much higher cost than our method, see Fig 2. 
The second part of Example 2 is to compare 
the composite rule described in Section 5
with the basic quadrature (\ref{es:s12}) of Section 4
when the length of the interval of integration, $[a,b]$, becomes 
increasingly large.
Here $M$ denotes the number of subintervals in 
$[0, T]$ and $Mn$ stands for the total number of support
points in the given interval $[0, T].$ 
\begin{figure}[ht]
\begin{center}
\vspace{0in}
\leavevmode
\hspace{0in}
\epsfxsize =3in
\epsfbox{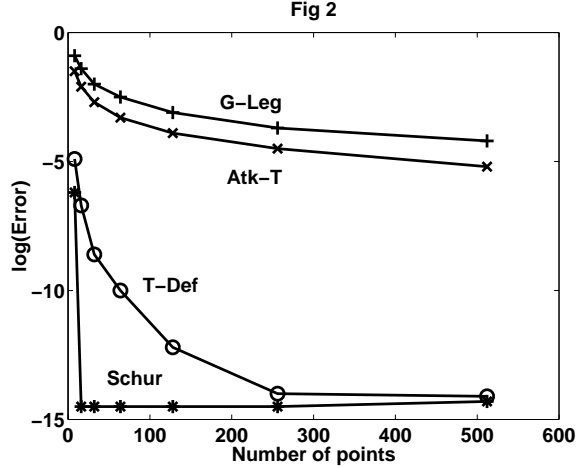}
\end{center}
\vspace{0in}
\hspace{0in}
\caption{ \label{figure2} Comparison of numerical solutions of Example 2.   }
\end{figure}
\vspace{0.2in}
\begin{center}
Table 1 ( \ $T=200\pi,$ \ \ $\lambda=-\frac{4}{\pi}$ \ ) \\ 
\begin{tabular}{|c|c|c|c|c|c|c|} \hline
$M$ & 1 & 1 & 1 & 2 & 4 & 8     \\ \hline
$Mn$ & $256$ & $512$ & $1024$ & $256$ & $512$ & $1024$  \\ \hline 
Error & $2.4e+01$ & $3.0e-02$ & CPUtime limit exceed & $1.2e+01$ & $7.8e-02$ & $2.2e-11$ \\ \hline 
\end{tabular}
\end{center}    
\vspace{0.2in} 
Without partitioning, i.e. with $M=1$,  we increase the number 
of support points
from $n=128$ to $n=1024.$ \ For $n=512$ the accuracy is of order
${\it {O}}(10^{-2}),$ \ but for $n=1024$ \  the CPU time limit is exceeded.
When the interval is partitioned into $8$ subintervals and $n=128$ 
i.e., the total number of points is $1024,$ 
the accuracy now is of order ${\it {O}}(10^{-11}).$ 
\\
{\bf{Example}} 3. 
\begin{eqnarray*}
x(t)+\int^{1}_{-1}k(t,s)x(s)ds=y(t), \ \ \ -1 \leq t \leq 1
\end{eqnarray*} 
where $y(t)=1-t^{2}+\frac{1}{1-t^{2}}(\arctan(t)-
 \arctan(-1))-\frac{1}{(1+t)(1+t^{2})},$
and \\
\[
k(t,s)  = 
\left\{ \begin{array}{cc}
\frac{1}{(1-t^{2})(1-s^{4})} & {\rm if} \ \ -1 \leq s \leq t \\
\frac{-1}{(1-t^{4})(1-s^{2})} & {\rm if} \ \ \ \  t < s \leq 1.
\end{array} 
\right.
\]
The analytical solution is $x(t)=1-t^2$.
Since this kernel has singularities along the boundaries of the
square $[-1,1]\times [-1,1]$ 
methods based on the Trapezium rule 
are not applicable. Therefore we compare  our algorithms of Section 1
and Section 3
with the $Nystr{\ddot{o}}m$-Gauss-Legendre discretization only.
The algorithm of Section 1 
shows the same  accuracy of numerical solution as 
the Gauss-Legendre quadrature. 
The method of Section 3  gives ${\it{O}}(10^{-13})$ accuracy with $n=32$
points, whereas  $Nystr{\ddot{o}}m$-Gauss-Legendre quadrature
gives  ${\it{O}}(10^{-3})$ with $n=256$ points.
\begin{figure}[ht]
\begin{center}
\vspace{0in}
\leavevmode
\hspace{0in}
\epsfxsize =3in
\epsfbox{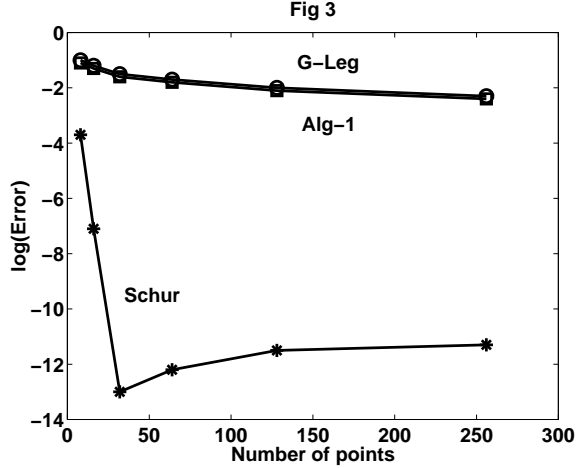}
\end{center}
\vspace{0in}
\hspace{0in}
\caption{ \label{figure3} Comparison of numerical solutions of Example 3.   }
\end{figure}
\vspace{0.2in}
\\
{\bf{Example}} 4.
\begin{eqnarray*}
x(t)+\int^{1}_{-1}k(t,s)x(s)ds=y(t), \ \ \ -1 \leq t \leq 1
\end{eqnarray*} 
where $y(t)=2(1-t^{2}+2t^{3})+(1+2t^{4})ln(t^2+t^4)
-ln(1+t^2)-2t^{4}ln(1+t^4),$ \
and \\ 
\[
k(t,s)=
\left\{ \begin{array}{cc}
1/(t^{2}+s^{4}) & {\rm if} \ \  -1 \leq s \leq t\\
1/(s^{2}+t^{4}) & {\rm if} \ \   t < s \leq 1.
\end{array}
\right.
\]
The analytical solution is $x(t)=4t^3$.
The kernel $k(t,s)$ has a singularity at $(0,0)$.  
Also $y(t)$ is singular at $t = 0$. 
 Since Chebyshev points $cos(\frac{(2i-1)}{2N}\pi)$, \
$i=1,2,...,n,$ are clustered towards the end points of interval $[-1, 1],$
discretization formula (\ref{es:s12}) does not contain  
sufficient values of kernel
near $(t,s)=(0,0)$. Therefore we partition $[-1,1]$ into  $[-1,0]$
and $[0,1].$  The choice of $n=256$  Chebyshev points in each subinterval 
with the total of $n=512$ points
gives ${\it O}(10^{-11})$
accuracy.
For comparison the best accuracy of the Gauss-Legendre quadrature
without partitions and with $n=512$ support points is
${\it{O}}(10^{-4})$ while the best accuracy of the algorithm
of Section 2 without partitions is of ${\it O}(10^{-6}),$ see Fig 4. 
\begin{figure}[ht]
\begin{center}
\vspace{0in}
\leavevmode
\hspace{0in}
\epsfxsize =3in
\epsfbox{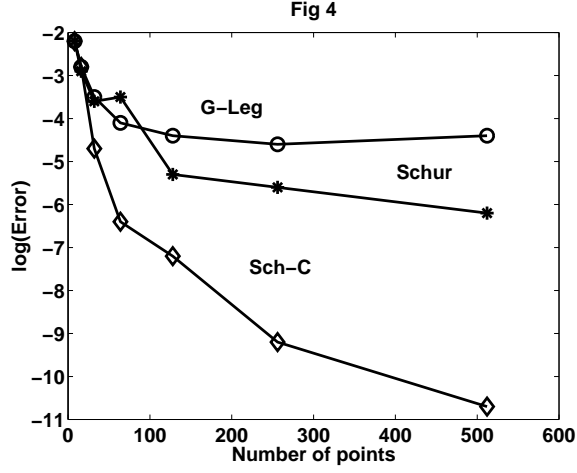}
\end{center}
\vspace{0in}
\hspace{0in}
\caption{ \label{figure4} Comparison of numerical solutions of Example 4.  }
\end{figure}
\vspace{0.2in}
\section{Application to Non-Local Schroedinger Equations}

In this section we demonstrate that the developed numerical 
technique is also applicable to problems other than integral equations,
for example, to integro-differential equations. We chose here 
the radial Schroedinger equation which models
the quantum mechanical interaction between particles represented by
spherically 
symmetric potentials.
These potentials are usually local, i.e., they depend only on the 
distance between the two particles, in which case the
equation is a differential equation which is
routinely solved in computational physics. 
However, if there are more than 
two particles present, then the potentials can become non-local and the 
differential Schroedinger equation becomes an 
integro-differential equation for the wave function $\psi,$   
\begin{equation}
\label{eq:pr1}
\frac{d^{2}\psi(r)}{dr^{2}}+{\kappa}^{2}\psi(r)=
\int^{T}_{0}v(r,r')\psi(r')dr',
\end{equation}
which is defined for $0 < r < \infty$, satisfies  the condition $\psi(0)=0,$ 
and is bounded at infinity.
It is assumed that $v(r,r')$ is negligible for 
$r>T$  or $r'>T,$ 
see e.g. \cite{mott}.
Because it is numerically more difficult to solve the Schroedinger equation
in the presence of a 
nonlocal potential, the latter is customarily replaced by an approximate
local equivalent potential.
There is, however, a renewed
interest in the nonlocal equations, and a  significant number of papers 
on this subject appeared in the past few years, (our database search
returned over 50 related publications).

Using the technique of \cite{reo} it easy to show that 
(\ref{eq:pr1}) is equivalent to the following integral equation,
\begin{eqnarray*}
\psi(r) & + & \frac{cos(\kappa r)}{\kappa}\int^{r}_{0}sin(\kappa r')
\int^{T}_{0}v(r',p)\psi(p)dpdr' \\
& + & \frac{sin(\kappa r)}{\kappa}\int^{T}_{r}cos(\kappa r')
\int^{T}_{0}v(r',p)\psi(p)dpdr'=sin({\kappa}r).
\end{eqnarray*} 
or
\begin{equation}
\label{eq:q1}
\psi(r) + \frac{cos(\kappa r)}{\kappa}\int^{T}_{0}
k_{1}(r,r')\psi(r')dr'+ 
\frac{sin(\kappa r)}{\kappa}\int^{T}_{0}k_{2}(r,r')\psi(r')dr'
\end{equation}
$$=sin(\kappa r),$$
where
\begin{eqnarray*}
k_{1}(r,r')=\int^{r}_{0}sin(\kappa p)v(p,r')dp, \ \ \
k_{2}(r,r')=\int^{T}_{r}cos(\kappa p)v(p,r')dp.
\end{eqnarray*}
We consider now the case when $v(p,r')$ is $(p_{1},p_{2})$-semismooth,
such that
\[
v(p,r')=
\left \{ \begin{array}{cc}
v_{1}(p,r') & {\rm if} \ \  0 \leq p \leq r'  \\
v_{2}(p,r') & {\rm if} \ \   r' \leq p \leq T .
\end{array}
\right.
\]
In order to use the method which we developed in previous sections, 
we rewrite 
equation (\ref{eq:q1}) as follows,
\begin{eqnarray*}
\psi(r) + \frac{c(r)}{\kappa}\int^{r}_{0}k_{1}(r,r')\psi(r')dr'+
\frac{c(r)}{\kappa}\int^{T}_{r}k_{1}(r,r')\psi(r')dr'
\end{eqnarray*}
\begin{equation}
\label{eq:q2}
+\frac{s(r)}{\kappa}\int^{r}_{0}k_{2}(r,r')\psi(r')dr'+
\frac{s(r)}{\kappa}\int^{T}_{r}k_{2}(r,r')\psi(r')dr' = s(r),
\end{equation}
where for notational convenience we abbreviate,
\ $c(r)=cos(\kappa r), s(r)=sin(\kappa r).$ 
We have
\[
k_{1}(r,r')=
\left \{ \begin{array}{cc}
k_{11}(r,r') & {\rm if} \ \  0 \leq r' \leq r  \\
k_{12}(r,r') & {\rm if} \ \  0 \leq r \leq r' , 
\end{array}
\right.
\]
and
\[
k_{2}(r,r')=
\left \{ \begin{array}{cc}
k_{21}(r,r') & {\rm if} \ \   r' \leq r \leq T  \\
k_{22}(r,r') & {\rm if} \ \   r \leq r' \leq T , 
\end{array}
\right.
\]
where
\begin{eqnarray*}
k_{11}(r,r') & = & \int^{r'}_{0}s(p)v_{1}(p,r')dp
+\int^{r}_{r'}s(p)v_{2}(p,r')dp, \\
& = & \int^{r'}_{0}s(p)v_{1}(p,r')dp
+\int^{r}_{0}s(p)v_{2}(p,r')dp-
\int^{r'}_{0}s(p)v_{2}(p,r')dp, \\
k_{12}(r,r') & = & \int^{r}_{0}s(p)v_{1}(p,r')dp \\
k_{21}(r,r') & = & \int^{T}_{r}c(p)v_{2}(p,r')dp \\
k_{22}(r,r') & = & \int^{r'}_{r}c(p)v_{1}(p,r')dp
+\int^{T}_{r'}c(p)v_{2}(p,r')dp,\\
& = & \int^{T}_{r}c(p)v_{1}(p,r')dp-
\int^{T}_{r'}c(p)v_{1}(p,r')dp
+\int^{T}_{r'}c(p)v_{2}(p,r')dp.
\end{eqnarray*}
Thus,
\begin{eqnarray*}
\psi(r) + \frac{c(r)}{\kappa}\int^{r}_{0}k_{11}(r,r')\psi(r')dr'+
\frac{c(r)}{\kappa}\int^{T}_{r}k_{12}(r,r')\psi(r')dr'
\end{eqnarray*}
\begin{equation}
\label{eq:q3}
+\frac{s(r)}{\kappa}\int^{r}_{0}k_{21}(r,r')\psi(r')dr'+
\frac{s(r)}{\kappa}\int^{T}_{r}k_{22}(r,r')\psi(r')dr' = s(r),
\end{equation}
Applying our quadrature to this equation
we get,
\begin{equation}
\label{eq:q5}
[ \ {\bf{I}}+\frac{T}{2\kappa}{\bf{D}}_{c}({\bf{W}} \circ {\bf{K}}_{11}
+{\bf{V}}\circ {\bf{K}}_{12})+\frac{T}{2\kappa}{\bf{D}}_{s}
({\bf{W}}\circ {\bf{K}}_{21}+{\bf{V}}\circ {\bf{K}}_{22}) \ ]
{\bar{\bf{\psi}}}
=
{\bar{\bf{s}}},
\end{equation}
where in more detail,
\begin{eqnarray*}
{\bar{\bf{\psi}}} & = & [\psi(t_{0}),\psi(t_{1}),...,\psi(t_{n})]^{T},\\
{\bf{D}}_{c} & = & diag(cos(\kappa t_{0}),cos(\kappa t_{1}),...,
cos(\kappa t_{n})),\\
{\bf{D}}_{s} & = & diag(sin(\kappa t_{0}),sin(\kappa t_{1}),...,
sin(\kappa t_{n})),\\
{\bf{D}}_{\sigma} & = & diag([1,...,1]{\bf{S}}_{L}{\bf{C}}^{-1}), \\
{\bar{\bf{s}}} & = & [sin(\kappa t_{0}),sin(\kappa t_{1}),...,
sin(\kappa t_{n})]^{T}, \\
{\bf{W}} & = & {\bf{C}}{\bf{S}}_{L}{\bf{C}}^{-1}, \ \ \
{\bf{V}} = {\bf{C}}{\bf{S}}_{R}{\bf{C}}^{-1}, \\
{\bf{K}}_{11} & = & (  k_{11} (t_{i},t_{j}) )^{n}_{i,j=0} \\
 & = & \frac{T}{2}[ diag (  {\bf{W}}{\bf{D}}_{s}
({\bf{V}}_{1}-{\bf{V}}_{2}) ) + 
 ( {\bf{W}}{\bf{D}}_{s}{\bf{V}}_{2}  ) ] \\
{\bf{K}}_{12} & = & (  k_{12} (t_{i},t_{j}) )^{n}_{i,j=0} \\ 
& = & \frac{T}{2} ({\bf{W}}{\bf{D}}_{s}{\bf{V}}_{1}),\\
{\bf{K}}_{21} & = & (  k_{21} (t_{i},t_{j}) )^{n}_{i,j=0}, \\  
& = & \frac{T}{2} ({\bf{V}}{\bf{D}}_{c}{\bf{V}}_{2}) \\
{\bf{K}}_{22} & = & (  k_{22} (t_{i},t_{j}) )^{n}_{i,j=0} \\ 
& = & \frac{T}{2}[  (  {\bf{V}}{\bf{D}}_{c}{\bf{V}}_{1}  ) + 
diag( {\bf{V}}{\bf{D}}_{c}({\bf{V}}_{2}-{\bf{V}}_{1})  )].  
\end{eqnarray*}
We illustrate now the obtained discretization with two examples. In the first
example we use a prototype of the Yukawa potential, (e.g. \cite{landa}, 23.c),
which is simplified to a degree such that an analytic solution can be found.
In our terminology this potential is semiseparable. 
We note once more that the case of this semi-separable potential  could
be treated more easily by the 
techniques already presented in 
\cite{reo}, and we use it here only because the comparison with the analytic
solution is possible. \\
\noindent
{\bf {Example 1.}}
Let       
{\begin{eqnarray*}
v(p,r')=
\left \{ \begin{array}{cc}
\lambda e^{p-r'} & {\rm if} \ \  0 \leq p \leq r'  \\
\lambda e^{r'-p} & {\rm if} \ \   r' \leq p \leq T . 
\end{array}
\right.
\end{eqnarray*} It is easy to see that if $\psi(r)=e^{-r},$ then
the right-hand side has the form,
$$y(r)=(1-\frac{3{\lambda}{\kappa}}{4})e^{-r}
+\frac{3{\lambda}{\kappa}}{4}cos(r)-\frac{{\lambda}{\kappa}}{2}re^{-r}.$$
By comparing the analytical solution given above with the
numerical solution of (\ref{eq:q5}) at the discretization points,
we get the following relative errors
in the case of ${\lambda}=0.1, \ \ {\kappa}=1$ \ \ and \ $ T=20.$
\begin{center}
\begin{tabular}{|c|c|c|c|c|c|} \hline
$n$ & 16 & 32 & 64 & 128 & 256 \\ \hline
$Error$ & $1.2e+01$ & $3.4e-07$ & $8.1e-09$ & $3.4e-09$ & $6.0e-09$ \\ \hline
\end{tabular}
\end{center} 

In the second example we consider a more interesting case for which the 
techniques of \cite{reo} are not applicable. This time the non-locality
is a prototype of the optical model Perey-Buck potential, 
(e.g. \cite{fesh}, Ch. V.2).
 In our terminology this potential is semi-smooth,
but not semiseparable. \\
\noindent
{\bf{Example 2.}} Let
\[
v(p,r')=\frac{{\lambda}e^{-\frac{|r'-p|}{A}}}{1+e^{-\frac{|r'-p|}{A}}}
       =
\left \{ \begin{array}{cc}
\frac{\lambda e^{\frac{p-r'}{A}}}{1+e^{\frac{p-r'}{A}}} & {\rm if} \ \  0 \leq p \leq r'  \\
\frac{\lambda e^{\frac{r'-p}{A}}}{1+e^{\frac{r'-p}{A}}} & {\rm if} \ \   r' \leq p \leq T .
\end{array}
\right.
\]
Solving (\ref {eq:q5}) 
at $n$ shifted Chebyshev support points $t^{(n)}_{i}, i=1,...,n,$  and $2n$
 points $s^{(2n)}_{i}, i=1,...,2n,$  
we obtain  numerical solutions ${\psi}^{(n)} (r)$ and ${\psi}^{(2n)} (r)$,
respectively.

To get the values of  ${\psi}^{(2n)}(r)$ at $t^{(n)}_{i},$
we follow the procedure described in the
beginning of Section 6.
The error $e_{n}$ is obtained by comparison of the solutions
${\psi}^{(n)}$ \ and ${\psi}^{(2n)}$ \ as follows,
$$e_{n}=||{\psi}^{(2n)}(t^{(n)}_{i})-
{\psi}^{(n)}(t^{(n)}_{i})||_{\infty}/||{\psi}^{(2n)}(t^{(n)}_{i})||
_{\infty}.$$ 
Here we take $\lambda=0.1$,\ $\kappa=1,$ \ $A=100,$ \ and \ $T=20.$
\begin{center}
\begin{tabular}{|c|c|c|c|c|c|c|} \hline
$n$ & 8 & 16 & 32 & 64 & 128 & 256  \\ \hline
$e_{n}$ & $1.0e-0$ & $1.2e-03$ & $1.6e-09$ & $7.7e-15$ 
& $1.6e-14$ & $4.8e-14$ \\ \hline
\end{tabular}
\end{center} 
We see that for this choice of $\lambda$ the matrix is well conditioned 
and the double precision accuracy is obtained with 64 points.

\section{Summary and Conclusions}

In this paper, which is one of a sequence treating integral equations,
we describe a new  
method for solving integral equations for the case when kernel
can be discontinuous along the main diagonal. It has the following
advantages for a large class of such kernels: 
(i) for semismooth kernels it gives a much higher accuracy than
it was ever possible with standard Gauss type quadrature rules;
(ii) it is of comparable accuracy with Gauss type quadratures for 
smooth kernels; (iii) it exploits additional structure of the kernel
such as a low semi-rank, or a displacement structure, for example,
to allow for reduced complexity algorithms for the discretized equations.
(iv) the numerical examples provided in the present study illustrate
increased accuracy of our method compared to other more
conventional methods.

Our method is also applicable to other problems, such as the computation of 
eigenvalues and eigenfunctions of integral and differential operators and
solution of integro-differential equations.

Our method may find applications in quantum mechanical atomic and
nuclear physics problems, where the requirement of indistinguishability
of the electrons  leads to non localities in the potential contained
in the Schroedinger equation due to the presence of exchange
terms. These, in turn, lead to integro-differential equations which
are usually solved by iterative finite difference methods,
or by orthogonal function expansion methods. We plan to
compare our new method with some of the existing methods in future
investigations on more realistic examples.

\end{document}